\documentclass[11pt]{article}
\usepackage[utf8]{inputenc}
\usepackage[english]{babel}
\usepackage{amssymb,amsmath}
\usepackage[usenames, dvipsnames]{color}
\usepackage{amsthm}
\usepackage{cite}
\usepackage{hyperref}
\let\OLDthebibliography\thebibliography
\renewcommand\thebibliography[1]{
  \OLDthebibliography{#1}
  \setlength{\parskip}{0pt}
  \setlength{\itemsep}{0pt plus 0.3ex}
}
\title{Explicit solutions of certain orientable quadratic equations in free
groups}
\author{D.~Gon\c{c}alves\footnote{The author is partially supported by FAPESP-Funda\c c\~ao 
de Amparo a Pesquisa do Estado de S\~ao Paulo, Projeto Tem\'atico 
Topologia Alg\'ebrica, Geom\'etrica 2016/24707-4.}, T.~Nasybullov\footnote{The author is supported by the Research Foundation -- Flanders (FWO): postdoctoral grant  12G0317N and travel grant  V418018 N.}}
\date{}

\newtheorem{theorem}{{\scshape Theorem}}[section]
\newtheorem{lemma}[theorem]{{\scshape Lemma	}}
\newtheorem{corollary}[theorem]{{\scshape Corollary}}
\newtheorem{proposition}[theorem]{{\scshape Proposition}}
\newtheorem{remark}[theorem]{{\scshape Remark}}
\begin{document}
\maketitle 
\begin{abstract}
For $g\geq1$ denote by $F_{2g}=\langle x_1, y_1,\dots,x_g,y_g\rangle$ the free group on $2g$ generators and let $B_g=[x_1,y_1]\dots[x_g,y_g]$.
For $l,c\geq 1$ and elements $w_1,\dots,w_l\in F_{2g}$ we study 
orientable  quadratic equations of the form $[u_1,v_1]\dots[u_h,v_h]=(B_g^{w_1})^c(B_g^{w_2})^c\dots(B_g^{w_l})^c$ with unknowns $u_1,v_1,\dots,u_h,v_h$ and provide explicit solutions for them for the minimal possible number $h$.

In the particular case when  $g=1$, $w_i=y_1^{i-1}$ for $i=1,\dots,l$ and $h$ the minimal number which satisfies $h \geq l(c-1)/2+1$  we provide two types of solutions depending on the image of the subgroup $H=\langle u_1,v_1,\dots,u_h,v_h\rangle$ generated by the solution under the natural homomorphism $p:F_2\to F_2/[F_2,F_2]$: the first solution, which is  called  a primitive  solution, satisfies $p(H)=F_2/[F_2,F_2]$, the second solution satisfies 
$p(H) = \big\langle p(x_1),p(y_1^l)\big\rangle$. 

We also provide an explicit solution of the equation $[u_1,v_1]\dots[u_k,
v_k] = \big(B_1\big)^{k+l} \big({B_1}^{y}\big)^{k-l}$ for $k>l\geq0$ in $F_2$,  and prove that if $l\neq0$, then every solution of this equation is primitive.

As a geometrical consequence, for every solution we obtain a map $f:S_h\to T$ from the orientable surface $S_h$ of genus $h$ to the torus $T=S_1$ which has the minimal number of roots among all maps from the homotopy class of $f$. Depending on the number $|p(F_2):p(H)|$ such maps have fundamentally different geometric properties: in some cases they satisfy the Wecken property and in other cases not.

~
\end{abstract}
\section{Introduction and preliminaries}\label{sec1}
Let $G$ be a group and $S$ be a symmetric subset of $G$, i.~e. a subset such that $1\notin S$ and $S=S^{-1}$. For an element $a$ from $\langle S\rangle$ denote by $l_S(a)$ the minimal number $k$ such that $a$ is a product of $k$ elements from $S$. The number $l_S(a)$ is called the length of $a$ with respect to $S$. The numbers $l_S(a)$ and especially the value ${\rm sup}\big(l_S(a)~|~a\in\langle S\rangle\big)$ for different sets $S$ in different groups $G$ has been studied by various authors (see, for example,
\cite{Seg} and references therein). 

If $S$ is the set of all nontrivial commutators in $G$, then $l_S(a)$ is called the commutator length (or the genus) of an element $a\in [G,G]$. The problem of determining the commutator length of an element $a\in [G,G]$ corresponds to the problem of finding the minimal $h$ for which  the equation
\begin{equation}\label{eqgen}
[u_1,v_1] \dots [u_h, v_h]= a 
\end{equation}
with unknowns $u_1,v_1,\dots,u_h,v_h$ admits a solution in $G$. Such equation is called an orientable quadratic equation. The word ``quadratic'' means that every variable in the left side of the equation appears exactly twice. The word ``orientable'' means that every unknown variable $x$ appears once in the form $x$ and once in the form $x^{-1}$.  If there exists a variable $x$ in a quadratic equation which appears twice with exponent $1$, then this equation is called non-orientable.
 
 The problem of finding a solution of equation (\ref{eqgen}) in the free group is closely related with coincidence theory of maps between  orientable surfaces, which includes the case of the study of roots.
  See, for example, 
   \cite[Fundamental Lemma 1.2]
 {Gon} and  \cite{BogGonKudZie2, BogGonZie, GonZie}. 

Many works concerning the problem of finding solutions of quadratic equations and especially of equation (\ref{eqgen}) in different groups have been done from both algebraic and geometric points of view. For many years great attention was paid to equations in free groups \cite{Cul, Gon,  BogGonKudZie2, BogGonZie, GonZie,  Hme1, Hme2, Hme3, KudWeiZie, Lyn1, Lyn2, Mak, Ste, Vdo, Wic1}.
 The particular case of equation (\ref{eqgen}) when $h=1$ was studied in \cite{Wic1} (see also \cite{GonKudZie3}). For a given quadratic equation with any number of unknown variables in any free group with the right-hand side an arbitrary element an algorithm for solving the problem of the existence of a solution was given by Culler \cite{Cul} using the surface method and generalizing the result of Wicks \cite{Wic1}. Based on different techniques, the problem has been studied by the first named author with coauthors \cite{GonKudZie1, GonKudZie2, GonKudZie3} for parametric families of quadratic equations arising from continuous maps between closed surfaces.   

The question about the existence of a solution of equation (\ref{eqgen}) can be solved in many cases. However, the majority of results are either about non-existence of a solution, or about existence only and they do not provide an algorithm how to find an explicit solution. The following result from \cite[Proposition 4.2]{GonZie} gives one simple necessary condition     for solvability of equation~(\ref{eqgen}) with the right-hand side of the special form motivated by geometry.

\begin{proposition}\label{nonex} Let $w_1,\dots,w_l$ be distinct elements of the free group $F_{2g}=\langle x_1, y_1, \dots,x_g,y_g\rangle$ and let $c_1,\dots,c_l$ be integers which are all positive or all negative. Denote by $B_g=[x_1,y_1]\dots [x_g,y_g]$. If the equation
\begin{equation}\label{eqpart}
[u_1,v_1] \dots [u_h,v_h]
 = \Big(B_g^{w_1}\Big)^{c_1} \Big(B_g^{w_2}\Big)^{c_2} \dots
\Big(B_g^{w_l}\Big)^{c_l}
\end{equation}
with unknowns $u_1,v_1,\dots,u_h,v_h$ is solvable in $F_{2g}$, then $(|c_1|+\dots+|c_l|)(2g-1)\leq 2h-2+l$.
\end{proposition}
\noindent Here and throughout the paper for elements $a,b$ we denote by $a^b=bab^{-1}$ the conjugate of $a$ by $b$, and by $[a,b]=aba^{-1}b^{-1}$ the commutator of elements $a,b$. 

Orientable quadratic equations with the right-hand side as in (\ref{eqpart}) (not necessary for $c_1,\dots,c_l$ all positive or all negative) are of special interest in geometry. If $f: S_h\to S_g$ is a continuous map between orientable surfaces, then it induces a map $f_{\#}: \pi_1(S_h)\to \pi_1(S_g)$ between fundamental groups. Denoting by $\pi_1(S_h)=\langle x_1,y_1,\dots,x_h,y_h~|~[x_1,y_1]\dots[x_h,y_h]=1\rangle$ and $f_{\#}(x_i)=u_i$, $f_{\#}(y_i)=v_i$ we must have $[u_1,v_1]\dots[u_h,v_h]=1$ in $\pi_1(S_g)$, i.~e. $[u_1,v_1]\dots[u_h,v_h]$ must be expressible as a right-side of~(\ref{eqpart}). So, there is a strong connection between maps between orientable surfaces and orientable quadratic equations with the right hand side as in (\ref{eqpart}).

Note that the result of Proposition \ref{nonex} says that if $(|c_1|+\dots+|c_l|)(2g-1)> 2h-2+l$, then equation (\ref{eqpart}) does not have a solution independently of  elements  $w_1,...,w_{l}$. However there is no guarantee that a solution exists for $(|c_1|+\dots+|c_l|)(2g-1)\leq 2h-2+l$. Moreover  a solution can exist for some elements $w_1,\dots,w_{l}$ but not for others. For example, if $g=1$, $l=2$, $c_1=c_2=1$, then $h=1$. Using Wicks criterion \cite{Wic1} it is easy to show that the equation $[u,v] = B_1 B_1^{y_1^2}$ has no solutions in $F_2=\langle x_1,y_1\rangle$. However, the equation $[u,v]=B_1 B_1^{y_1}$ has a solution $u=x_1$, $v=y_1^2$.
So, it is reasonable to ask for which integers $c_1,\dots,c_{l}$ and elements $w_1,\dots,w_{l}$   equation $(\ref{eqpart})$  has a solution, and when it has, provide this solution. Not much is known about this problem.

In the present work we consider equation (\ref{eqpart}) in the free group $F_{2g}$  with right parts of special forms  and our goal is  to provide explicit solutions for them. In turn, for $g=1$ this will provide existence of maps from the orientable surface of genus $h$ into the torus which have some features about root theory. For some cases we will find two types of solutions depending on the index of the image of the subgroup $H=\langle u_1,v_1,\dots,u_h,v_h\rangle$ generated by the solution under the natural homomorphism $p:F_{2g}\to \pi_1(S_g)=F_{2g}/\langle B_g \rangle^{F_{2g}}$ in $\pi_1(S_g)$. In order to explain the importance of the value $|\pi_1(S_g):p(H)|$ let us recall some facts from Nielsen root theory.

Let $f:M_1\to M_2$ be a continuous map between closed manifolds $M_1$, $M_2$ and let $c\in M_2$. Every element from $f^{-1}(c)$ is called a root. The minimal number of roots in  the homotopy class of a map $f$ is the number $MR[f]={\rm min}_{g\simeq f}\big(|g^{-1}(c)|\big)$, where $\simeq$ denotes the homotopy equivalence. This number does not depend on $c$. Two roots $x,y\in M_1$ are said to belong to the same Nielsen root class if there exists a path $\gamma$ in $M_1$ connecting $x,y$ such that $f(\gamma)$ is contractible. For a map $f$ between two manifolds of the same dimension an index for a Nielsen root class is defined in \cite{Ki}. A Nielsen root class is called essential if its index is not equal to zero. The indices of all essential Nielsen root classes coincide. The Nielsen root number $NR[f]$ is the number of essential Nielsen root classes, this number is always finite and it satisfies the inequality $NR[f]\leq MR[f]$.  If $NR[f]=MR[f]$, then $f$ is said to possess the Wecken property. The map $f$ induces the map $f_{\#}:\pi_1(M_1)\to \pi_1(M_2)$ between fundamental groups. Denote by $l(f)=|\pi_1(M_2):f_{\#}(\pi_1(M_1))|$ if $|\pi_1(M_2):f_{\#}(\pi_1(M_1))|$ is finite, and $l(f)=0$ otherwise. If $M_1=S_h$, $M_2=S_g$ are closed orientable surfaces of genus $h,g$ respectively, then the map $f$ induces a homomorphism beween second homology groups $\mathbb{Z}=H_2(S_h)\to H_2(S_g)=\mathbb{Z}$. This map acts as a multiplication by some number $n$. This number is called the degree of $f$ and is denoted by ${\rm deg}(f)$. The absolute value $|{\rm deg}(f)|$ is denoted by $A(f)$. In \cite[Theorem 1.1]{BogGonKudZie2} it is proved that, if $A(f)\neq 0$, then 
\begin{align}\label{wecken}
MR[f]={\rm max}\Big(l(f), \chi(M_1)+(1-\chi(M_2))A(f)\Big)&&NR[f]=l(f),
\end{align} 
where $\chi$ denotes the Euler characteristic of the surface.

If  $u_1,v_1,\dots,u_h,v_h$ is some solution of equation (\ref{eqpart}), then one can construct a continuous map $f:S_h\to S_g$ which satisfies the following conditions: ${\rm deg}(f)=c_1+\dots+c_l$, $|f^{-1}(y)|=l$ for some point $y\in S_g$, the index of every Nielsen root class of $f$ is equal to $c_{i_1}+\dots+c_{i_k}$ for some indices $i_1,\dots,i_k$ and if $\pi_1(S_h)=\langle x_1,y_1,\dots,x_h,y_h~|~[x_1,y_1]\dots[x_h,y_h]=1\rangle$, then $f_{\#}(x_i)=u_i$, $f_{\#}(y_i)=v_i$. If  ${\rm deg}(f)=0$, then  $NR[f]=0$ and   $|p(F_{2g}):p(H)|$ can be either finite or infinite. If ${\rm deg}(f)\neq 0$,  then   $|p(F_{2g}):p(H)|$ is finite and   $NR[f]=|p(F_{2g}):p(H)|$, where $p:F_{2g}\to \pi_1(S_g)$ is a natural homomorphism and $H=\langle u_1,v_1,\dots,u_h,v_h\rangle$. See details about the construction of $f$ in \cite[Proposition 4.2]{GonZie}, here we are going to use only the properties of the constructed map.

The following result gives some information about the index of $p(H)$ in $\pi_1(S_g)$.
\begin{proposition}\label{crit} Let $p:F_{2g} \to \pi_1(S_g)$ be the homomorphism which sends the free generators of $F_{2g}$ to the canonical system of generators of the fundamental group $\pi_1(S_g)$ of an orientable surface $S_g$ of genus $g$. If $u_1,v_1,\dots,u_h,v_h$ is a solution of equation (\ref{eqpart}) with $c_1+c_2+...+c_{l}\neq 0$, then the index of $p\big(\langle u_1,v_1,\dots,u_h,v_h\rangle\big)$ in $\pi_1(S_g)$ is less than or equal to $l$.
\end{proposition}
\noindent\textbf{Proof.} Let $f:S_h \to S_g$ be the described before the proposition map constructed by the solution $u_1,v_1,\dots,u_h,v_h$. For some point $y\in S_g$ the number of elements in the preimage of $y$ under $f$ is $l$, therefore $MR[f]\leq l$.
Since ${\rm deg}(f)=c_1+\dots+c_l\neq 0$, we have $|p(F_{2g}):p(H)|=NR[f]\leq MR[f]\leq l$.\hfill$\Box$

In the present paper we will find explicit solutions for particular cases of equation (\ref{eqpart})  which are in some sense ``critical'' from the point of view of Proposition \ref{crit}: the first solution which is called a primitive solution satisfies the equality $p(H)=\pi_1(S_g)$ (the word ``primitive'' appears here naturally from the notion of primitives in free groups \cite{OsbZie}), and the second solution satisfies $|\pi_1(S_g):p(H)|=l$.

In Section \ref{moregen} we study the equation (\ref{eqpart}) for $c_1=c_2=\dots=c_l=c$ and prove that if for $c=1$ this equation has a solution which generates a subgroup $H_1$ of $F_{2g}$, then for every $c$ it has a solution which generates a subgroup $H_2$ such that $p(H_1)=p(H_2)$ (Theorem \ref{maingen}). In Section~\ref{var1} we consider the particular case of this equation in the free group $F_2=\langle x,y\rangle$ with the right part of the form 
$\big([x,y]\big)^c \big({[x,y]}^y\big)^c\dots \big({[x,y]}^{y^{l-1}}\big)^c$
 for integers $c, l\geq1$. We provide an explicit algebraic algorithm for finding the solution in a minimal subgroup (Corollary \ref{main2}), and an algorithm for finding a primitive solution (Theorem \ref{main1}) for such equation. In Section \ref{var2} we consider equation (\ref{eqpart}) in the free group $F_2=\langle x,y\rangle$ with the right part of the form $\big([x,y]\big)^{k+l} \big({[x,y]}^y\big)^{k-l}$ for integers $k>l\geq0$. We construct an explicit primitive solution for such equation (formulas (\ref{for3}), (\ref{for4})) 
 and prove that for $l\neq0$ every solution of such equation is primitive (Theorem \ref{main3}). Some geometrical consequences derived from this algebraic results are formulated in Corollaries~\ref{main02},~\ref{main01}.

{\bf Acknowledment:}  The first author would like to thank Prof. Richard  Weidmann for many helpful  discussions about the subject of the content of this work. 

\section{The right part has the form $\Big(B_g^{w_1}\Big)^c\Big(B_g^{w_2}\Big)^c\dots \Big(B_g^{w_l}\Big)^c$}\label{moregen}
The purpose of this section is to give an explicit solution for  equation (\ref{eqpart}) in the particular case when $c_1=c_2=\dots=c_l=c$
\begin{equation}\label{eeqq4}
[u_1,v_1]\dots[u_h,v_h] = \Big(B_g^{w_1}\Big)^c \Big({B_g}^{w_2}\Big)^c\dots \Big({B_g}^{w_l}\Big)^c
\end{equation}
\noindent for   the   minimal integer $h$ which satisfies the inequality
$h \geq l(c(2g-1)-1)/2+1$. We can assume that $c>0$ since if $c<0$, then denoting by $x_i^{\prime}=y_{g+1-i}$, $y_i^{\prime}=x_{g+1-i}$ for $i=1,\dots,g$ we have $F_{2g}=\langle x_1^{\prime}, y_1^{\prime}, \dots, x_g^{\prime}, y_g^{\prime}\rangle$, and in these generators equation (\ref{eeqq4}) has the same form where $c$ is changed by $-c$.  In the case when $h < l(c(2g-1)-1)/2+1$  by Proposition \ref{nonex} equation (\ref{eeqq4}) does not have solutions.

At first, we need the following simple lemma.
\begin{lemma}\label{tech} The word $w=a \xi_1 b \xi_2 c$ is a product of the commutator $[a \xi_1 a^{-1}, aba^{-1}]$ and the element $ab\xi_1\xi_2 c$. 
\end{lemma}
\noindent \textbf{Proof.} Straightforward calculation.
\hfill $\Box$

The main result of this section is the following theorem.
\begin{theorem}\label{maingen}Let $l,c,g\geq 1$ be integers, $h$ be the minimal integer which satisfies the inequality $h\geq l\big(c(2g-1)-1\big)/2+1$, $w_1,\dots,w_l$ be elements of the free group $F_{2g}=\langle x_1, y_1, \dots,x_g,y_g\rangle$ and $B_g=[x_1,y_1]\dots [x_g,y_g]$. If for $c=1$ the equation
\begin{equation}\label{eqpart2}
[u_1,v_1] \dots [u_h,v_h]
 = \Big(B_g^{w_1}\Big)^{c} \Big(B_g^{w_2}\Big)^{c} \dots
\Big(B_g^{w_l}\Big)^{c}
\end{equation}
has a solution  which generates the subgroup $H_1$ of $F_{2g}$, then for an arbitrary $c\geq 1$ it has a solution (explicitly constructed from the given solution for $c=1$)  which generates the subgroup $H_2$ of $F_{2g}$ such that if $p:F_{2g}\to\pi_1(S_g)$ is the natural homomorphism, then $p(H_1)=p(H_2)$.
\end{theorem}
\noindent \textbf{Proof.} We will construct the solution inductively on the variable $c$. The basis of induction $c=1$ is given as the condition of the theorem. Suppose that the statement is proved for $c=n$ and let us prove that it holds for $c=n+1$. We will consider two cases depending on the parity of $l$.

\textit{Case 1: $l$ is even.} Rewrite the right-hand side of equation (\ref{eqpart2}) for $c=n+1$ in the following form.
\begin{align}
\notag \big(B_g^{w_1}\big)^{n+1} \big(B_g^{w_2}\big)^{n+1} \dots
\big(B_g^{w_l}\big)^{n+1}&=\big(B_g^{w_1}\big)\Big(\big(B_g^{w_1}\big)^{n} \big(B_g^{w_2}\big)^{n} \dots
\big(B_g^{w_l}\big)^{n}\Big)\big(B_g^{w_1}\big)^{-1}\\
\notag&\cdot\big(B_g^{w_1}\big)\big(B_g^{w_l}\big)^{-n} \big(B_g^{w_{l-1}}\big)^{-n} \dots
\big(B_g^{w_2}\big)^{-n}\\
\label{genlev}&\cdot\big(B_g^{w_2}\big)^{n+1} \big(B_g^{w_3}\big)^{n+1} \dots
\big(B_g^{w_l}\big)^{n+1}
\end{align}

By the induction hypothesis there exist $u_1,v_1,\dots,u_h,v_h$ for $h=l\big(n(2g-1)-1\big)/2+1$ such that $[u_1,v_1] \dots [u_h,v_h]
 = \big(B_g^{w_1}\big)^n\big(B_g^{w_2}\big)^n \dots
\big(B_g^{w_l}\big)^n$ and $p\big(\langle u_1,v_1,\dots,u_h,v_h\rangle\big)=p(H_1)$. So, it is enough to prove that the product of two last lines of equation (\ref{genlev})
\begin{equation}\label{rewrgenlev}
\big(B_g^{w_1}\big)\big(B_g^{w_l}\big)^{-n} \big(B_g^{w_{l-1}}\big)^{-n} \dots
\big(B_g^{w_3}\big)^{-n}\big(B_g^{w_2}\big)\big(B_g^{w_3}\big)^{n+1} \dots
\big(B_g^{w_l}\big)^{n+1}
\end{equation}
is the product of $l\big((n+1)(2g-1)-1\big)/2+1-l\big(n(2g-1)-1\big)/2-1=l(2g-1)/2$ commutators of elements images of which under $p$ belong to $p(H_1)$. In equation (\ref{rewrgenlev}) denoting by $a=B_g^{w_1}$, $\xi_1=\big(B_g^{w_l}\big)^{-n}\big(B_g^{w_{l-1}}\big)^{-n}$, $b=\big(B_g^{w_{l-2}}\big)^{-n}\dots \big(B_g^{w_3}\big)^{-n}\big(B_g^{w_2}\big)\big(B_g^{w_3}\big)^{n+1}\dots \big(B_g^{w_{l-2}}\big)^{n+1}\big(B_g^{w_{l-1}}\big)$, $\xi_2=\big(B_g^{w_{l-1}}\big)^{n}\big(B_g^{w_l}\big)^{n}$, $c=\big(B_g^{w_l}\big)$ and applying Lemma \ref{tech} we conclude that expression (\ref{rewrgenlev}) is a product of the commutator of elements which belong to the kernel of $p$ times the element 
$$
\big(B_g^{w_1}\big)\big(B_g^{w_{l-2}}\big)^{-n} \big(B_g^{w_{l-3}}\big)^{-n} \dots
\big(B_g^{w_3}\big)^{-n}\big(B_g^{w_2}\big)\big(B_g^{w_3}\big)^{n+1} \dots
\big(B_g^{w_{l-2}}\big)^{n+1}\big(B_g^{w_{l-1}}\big)\big(B_g^{w_l}\big).
$$
Repeating this idea denoting by $\xi_1=\big(B_g^{w_{l-2}}\big)^{-n}\big(B_g^{w_{l-3}}\big)^{-n}$, $\xi_2=\big(B_g^{w_{l-3}}\big)^{n}\big(B_g^{w_{l-2}}\big)^{n}$, we conclude that expression (\ref{rewrgenlev}) is a product of two commutators of elements which belong to the kernel of $p$ times the element 
$$
\big(B_g^{w_1}\big)\big(B_g^{w_{l-4}}\big)^{-n} \big(B_g^{w_{l-5}}\big)^{-n} \dots
\big(B_g^{w_3}\big)^{-n}\big(B_g^{w_2}\big)\big(B_g^{w_3}\big)^{n+1} \dots
\big(B_g^{w_{l-4}}\big)^{n+1}\big(B_g^{w_{l-3}}\big)\dots\big(B_g^{w_l}\big).
$$
Repeating this procedure $(l-2)/2$ times we conclude that expression (\ref{rewrgenlev}) is the product of $(l-2)/2$ commutators of elements which belong to the kernel of $p$ times the element 
$B_g^{w_1}B_g^{w_2} \dots
B_g^{w_l}$ which (by the induction hypothesis for $c=1$) is the product of $l(g-1)+1$ commutators of elements images of which under $p$ belong to $p(H_1)$. Therefore expression
 (\ref{rewrgenlev}) is the product of $(l-2)/2+l(g-1)+1=l(2g-1)/2$ commutators of elements images of which under $p$ belong to $p(H_1)$.

\textit{Case 2: $l$ is odd.} Rewrite the right-hand side of equation (\ref{eqpart2}) for $c=n+1$ in the following form.
\begin{align}
\notag \big(B_g^{w_1}\big)^{n+1} &\big({B_g}^{w_2}\big)^{n+1}\dots \big({B}_g^{w_l}\big)^{n+1}=\\
\notag =&\Big(\big(B_g^{w_1}\big) \big({B_g}^{w_2}\big)\dots \big({B}_g^{w_l}\big)\Big)\\\notag \cdot&\Big(\big(B_g^{w_2}\big) \big({B_g}^{w_3}\big)\dots \big({B}_g^{w_l}\big)\Big)^{-1}\Big(\big(B_g^{w_1}\big) \big({B_g}^{w_2}\big)\dots \big({B}_g^{w_l}\big)\Big)\Big(\big(B_g^{w_2}\big) \big({B_g}^{w_3}\big)\dots \big({B}_g^{w_l}\big)\Big)\\
\notag \cdot&\Big(\big(B_g^{w_2}\big) \big({B_g}^{w_3}\big)\dots \big({B}_g^{w_l}\big)\Big)^{-2}\Big(\big(B_g^{w_1}\big)^{n-1}\big({B_g}^{w_2}\big)^{n-1}\dots \big({B}_g^{w_l}\big)^{n-1}\Big)\Big(\big(B_g^{w_2}\big) \big({B_g}^{w_3}\big)\dots \big({B}_g^{w_l}\big)\Big)^{2}\\
\notag\cdot& \big({B_g}^{w_l}\big)^{-1}\dots\big({B_g}^{w_2}\big)^{-1}\big({B_g}^{w_l}\big)^{-1}\dots\big({B_g}^{w_2}\big)^{-1}\\
\label{newc2}\cdot&\big({B_g}^{w_l}\big)^{-n+1}\dots \big({B_g}^{w_3}\big)^{-n+1}\big({B_g}^{w_2}\big)^{2}\big({B_g}^{w_3}\big)^{n+1}\dots \big({B_g}^{w_l}\big)^{n+1}
\end{align}
At first, we consider the particular case $n+1=2$. In this case equality (\ref{newc2}) implies
\begin{align}
\notag \big(B_g^{w_1}\big)^{2}& \big({B_g}^{w_2}\big)^{2}\dots \big({B}_g^{w_l}\big)^{2}=\\
\notag=&\Big(\big(B_g^{w_1}\big) \big({B_g}^{w_2}\big)\dots \big({B}_g^{w_l}\big)\Big)\\
\notag \cdot&\Big(\big(B_g^{w_2}\big) \big({B_g}^{w_3}\big)\dots \big({B}_g^{w_l}\big)\Big)^{-1}\Big(\big(B_g^{w_1}\big) \big({B_g}^{w_2}\big)\dots \big({B}_g^{w_l}\big)\Big)\Big(\big(B_g^{w_2}\big) \big({B_g}^{w_3}\big)\dots \big({B}_g^{w_l}\big)\Big)\\
\notag\cdot& \big({B_g}^{w_l}\big)^{-1}\dots\big({B_g}^{w_2}\big)^{-1}\big({B_g}^{w_l}\big)^{-1}\dots\big({B_g}^{w_2}\big)^{-1}\\
\label{newcc2}\cdot&\big({B_g}^{w_2}\big)^{2}\big({B_g}^{w_3}\big)^{2}\dots \big({B_g}^{w_l}\big)^{2}
\end{align}
So, if $n+1=2$, then since by induction hypothesis for $c=1$ the product of the first two lines of equation (\ref{newcc2}) is a product of $2l(g-1)+2$ commutators, it is enough to prove that the product of two last lines in~({\ref{newcc2}})
\begin{equation}\label{newcc22} \big({B_g}^{w_l}\big)^{-1}\dots\big({B_g}^{w_2}\big)^{-1}\big({B_g}^{w_l}\big)^{-1}\dots\big({B_g}^{w_3}\big)^{-1}\big({B_g}^{w_2}\big)\big({B_g}^{w_3}\big)^{2}\dots \big({B_g}^{w_l}\big)^{2}
\end{equation}
is a product of $\lceil l\big(2(2g-1)-1\big)/2\rceil+1-2l(g-1)-2=(l-1)/2$ commutators. Similarly to the first case (when $l$ is even) applying Lemma \ref{tech} to expression (\ref{newcc22}) for $\xi_1=\big(B^{w_3}\big)^{-1}\big(B^{w_2}\big)^{-1}$, $\xi_2=\big(B^{w_2}\big)\big(B^{w_3}\big)$ we conclude that (\ref{newcc22}) is a product of a commutator of elements which belong to $\langle B_g\rangle^{F_{2g}}$ times the same expression without elements $B^{w_2}, B^{w_3}$
$$\big({B_g}^{w_l}\big)^{-1}\dots\big({B_g}^{w_4}\big)^{-1}\big({B_g}^{w_l}\big)^{-1}\dots\big({B_g}^{w_5}\big)^{-1}\big({B_g}^{w_4}\big)\big({B_g}^{w_5}\big)^{2}\dots \big({B_g}^{w_l}\big)^{2}.
$$

Repeating this procedure $(l-1)/2$ times we conclude that expression (\ref{newcc22}) is a product of $(l-1)/2$ commutators of elements which belong to the kernel of $p$, i.~e. the case $n+1=2$ is proved.

For the general case $n>1$ by the induction hypothesis for $c=1$ the product of the first two lines of (\ref{newc2}) is the product of $2l(g-1)+2$ commutators, by the induction hypothesis for $c=n-1$ the third line of (\ref{newc2}) is the product of $\lceil l\big((n-1)(2g-1)-1\big)/2+1\rceil$ commutators of elements such that their images under $p$ generate $H_1$ (the case $n+1=2$ was necessary for making this inductive step). So, it is enough to prove that the product of two last lines in expression (\ref{newc2})
\begin{multline}\label{newc2rewr} \big({B_g}^{w_l}\big)^{-1}\dots\big({B_g}^{w_2}\big)^{-1}\big({B_g}^{w_l}\big)^{-1}\dots\big({B_g}^{w_2}\big)^{-1}\cdot\\
\cdot\big({B_g}^{w_l}\big)^{-n+1}\dots \big({B_g}^{w_3}\big)^{-n+1}\big({B_g}^{w_2}\big)^{2}\big({B_g}^{w_3}\big)^{n+1}\dots \big({B_g}^{w_l}\big)^{n+1}
\end{multline}
is a product of $\lceil l\big((n+1)(2g-1)-1\big)/2+1\rceil-\lceil l\big((n-1)(2g-1)-1\big)/2+1\rceil-2l(g-1)-2=l-2$
 commutators of elements images of which under $p$ belong to $p(H_1)$.

If we apply Lemma \ref{tech} twice to expression (\ref{newc2rewr}) for $\xi_1=\big(B_g^{w_3}\big)^{-1}\big(B_g^{w_2}\big)^{-1}$, $\xi_2=\big(B_g^{w_2}\big)\big(B_g^{w_3}\big)$, then we conclude that expression (\ref{newc2rewr}) is a product of two commutators (of elements which belong to the kernel of $p$) times expression (\ref{newc2rewr}) without elements $\big(B_g^{w_2}\big), \big(B_g^{w_3}\big)$. If we repeat this procedure $(l-3)/2$ times, we conclude that expression (\ref{newc2rewr}) is a product of $2(l-3)/2=l-3$ commutators times the expression
$$\big(B_g^{w_l}\big)^{-1}\big(B_g^{w_{l-1}}\big)^{-1}\big(B_g^{w_l}\big)^{-1}\big(B_g^{w_{l-1}}\big)^{-1}\big(B_g^{w_l}\big)^{-n+1}\big(B_g^{w_{l-1}}\big)^2\big(B_g^{w_l}\big)^{n+1}$$ 
which is equal to the commutator $\big(B_g^{w_l}\big)^{-1}\Big[\big(B_g^{w_{l-1}}\big)^{-1}\big(B_g^{w_l}\big)^{n-1},\big(B_g^{w_l}\big)^{-n}\big(B_g^{w_{l-1}}\big)^{-1}\Big]\big(B_g^{w_l}\big)$, i.~e. expression (\ref{newc2rewr}) is a product of $l-2$ commutators. 
\hfill $\Box$
\begin{remark}\label{rmain2} {\rm We can suppose that if $u_1,v_1,\dots,u_h,v_h$ is a solution of equation (\ref{eqpart2}) constructed in the proof of Theorem \ref{maingen} for an arbitrary $c$, then $u_1,v_1,\dots,u_{l(g-1)+1}, v_{l(g-1)+1}$ is the solution of equation (\ref{eqpart2}) for $c=1$.}
\end{remark}
\begin{remark} {\rm In particular case when $g=1$, $l=1$, $w_1=1$, $c=2h-1$ the explicit solution of equation (\ref{eqpart2}) is given in \cite[Proposition 4.6]{GonZie}.}
\end{remark}
The following  result, which is a consequence of the remark above, 
 is a corollary of Theorem~\ref{maingen}.
\begin{corollary} Let $G$ be a group, $a,b\in G$, and $n$ be an integer. Then $[a,b]^n$ can be expressed as the product of at most $\lceil (n+1)/2\rceil$ commutators.
\end{corollary}
\noindent\textbf{Proof.} The equation $[u,v]=[x,y]$ has a solution $u=x,v=y$ in $F_2=\langle x,y\rangle$. Therefore the equation $[u_1,v_1]\dots[u_h,v_h]=[x,y]^{n}$ has a solution in $F_2=\langle x,y\rangle$ for $h=\lceil (n+1)/2\rceil$. Acting on the equality $[u_1,v_1]\dots[u_h,v_h]=[x,y]^{n}$ by the homomorphism $\varphi:F_2\to G$ which is induced by $\varphi(x)=a$, $\varphi(y)=b$ we get the result.\hfill $\square$ 
\section{The right part has the form $\Big([x,y]\Big)^c \Big({[x,y]}^y\Big)^c\dots \Big({[x,y]}^{y^{l-1}}\Big)^c$}\label{var1}
In this section we will consider a particular case of equation (\ref{eqpart2}) in $F_2=\langle x,y\rangle$. In this case $p:F_2\to\pi_1(S_1)=F_2/[F_2,F_2]$ is the abelianization map. Denote by $B=B_1=[x,y]$. The following statement gives a stronger version of Proposition \ref{crit} for $c=g=h=1$.
\begin{proposition}\label{always}
Let $w_1,\dots,w_l$ be $l$ distinct elements from $F_2$. If $u,v$ is a solution of the equation 
$$[u,v]=B^{w_1}\dots B^{w_l},$$
and $p:F_2\to F_2/[F_2,F_2]$ is the natural homomorphism, then $|p(F_2):p\big(\langle u,v\rangle\big)|=l$.
\end{proposition}
\noindent \textbf{Proof.} By the solution $u,v$ we can construct a map $f:T\to T$ such that ${\rm deg}(f)=l\neq 0$. Since ${\rm deg}(f)\neq 0$, we have $NR[f]=|p(F_2):p\big(\langle u,v\rangle\big)|$. From the other side, since $f$ is a map from torus to torus, $NR[f]=|{\rm deg}(f)|$. 
This follows promptly from the main result in 
\cite{BrBPT} once one can identify the roots of $f$ with the fixed points of the map $g$ given by 
$g(x)=f(x)x$, using the multiplication of the torus.  \hfill$\square$

The purpose of this section is to give explicit solutions for the particular case of equation~(\ref{eqpart2})
\begin{equation}\label{eq1}
[u_1,v_1]\dots[u_h,v_h] = \Big([x,y]\Big)^c \Big({[x,y]}^y\Big)^c\dots \Big({[x,y]}^{y^{l-1}}\Big)^c
\end{equation}
\noindent in $F_2=\langle x,y\rangle$ for   the   minimal integer $h$ which satisfies the inequality
$h \geq l(c-1)/2+1$. If $h < l(c-1)/2+1$, then by Proposition \ref{nonex} equation (\ref{eq1}) does not have solutions.

In contrast to Proposition \ref{always}, for $c>1$ (and therefore $h>1$) the index $|p(F_2):p(H)|$, where $H=\langle u_1,v_1,\dots,u_h,v_h\rangle$, can be different. By Proposition \ref{crit} this index is at most $l$. We are going to introduce two types of solution of equation (\ref{eq1}): the first solution has the maximal possible index $|p(F_2):p(H)|=l$, and the second solution is primitive, i.~e. it has a minimal possible index $|p(F_2):p(H)|=1$.

Since for $c=1$ equation (\ref{eq1}) has a solution $u=x$, $v=y^l$, the case $|p(F_2):p(H)|=l$ is easy and it follows from  Theorem \ref{maingen} in the following way.
\begin{corollary}\label{main2}
Let $c,l \geq1$ be integers and $h$ be the minimal number which satisfies the inequality $h \geq l(c-1)/2+1$. Then the equation 
$$[u_1,v_1]\dots[u_h,v_h] = \Big([x,y]\Big)^c \Big({[x,y]}^y\Big)^c\dots \Big({[x,y]}^{y^{l-1}}\Big)^c$$
with unknowns $u_1,v_1,\dots,u_h, v_h$ has an explicit solution in $F_2$ given by recurrence in $c$ which satisfies the equality $p\Big(\big\langle u_1, v_1,\dots,u_h,v_h\big\rangle\Big) = \big\langle p(x),p(y^l)\big\rangle$, where $p:F_2\to F_2/[F_2,F_2]$ is the natural homomorphism.
\end{corollary}
Using equality (\ref{wecken}) and construction of the map $f$ (obtained from the solution) described in Section \ref{sec1} and in details in \cite[Proposition 4.2]{GonZie} we have the following corollary.  
\begin{corollary}\label{main02} Let $c>1$, $l\geq1$ be integers and $h$ be a minimal number which satisfies the inequality $h \geq l(c-1)/2+1$. Then there exists a map $f: S_h \to T$ with $A(f)=lc$,  $MR[f]=l$, $NR[f]=l$ 
and  each Nielsen root class has  index $c$. So, the Wecken property holds for $f$.
\end{corollary}

Now we are going to construct an explicit primitive solution of equation (\ref{eq1}). Results \cite{BogGonKudZie1, BogGonKudZie3} about primitive branching coverings give some evidence that such solutions might exist.
 At first, we consider one simple particular case when $l=2$, $c=2$.
\begin{lemma}\label{l2}
The equation 
$$[u_1,v_1][u_2,v_2]=[x,y]^2([x,y]^y)^2$$ 
has as solution  $u_1=x, v_1=y^3, u_2=y^3xy^{-2}x^{-1}y^{-1}xy^2x^{-1}y^{-3}, v_2=y^2x^2y^2x^{-1}y^{-3}$ which is primitive. 
\end{lemma}
\noindent\textbf{Proof.} Straightforward calculation.\hfill $\Box$

The general case follows.
\begin{theorem}\label{main1}Let $c>1$, $l\geq1$ be integers and $h$ be a minimal number which satisfies the inequality $h \geq l(c-1)/2+1$. Then the equation 
$$[u_1,v_1]\dots[u_h,v_h] = \Big([x,y]\Big)^c \Big({[x,y]}^y\Big)^c\dots \Big({[x,y]}^{y^{l-1}}\Big)^c$$
with unknowns $u_1,v_1,\dots,u_h,v_h$ has an explicit primitive solution in $F_2$ given by recurrence in~$l$. 
\end{theorem}
\noindent\textbf{Proof.} We will use induction on $l$. For the basis of induction we consider two cases $l=1$ and $l=2$. The result for $l=1$ follows from Corollary \ref{main2}. If $l=2$, then $h=c$. For $c=2$ the result follows from Lemma \ref{l2}. Suppose that we found a solution 
 $u_1,v_1,\dots,u_h,v_h$ for an integer $c$ such that $u_1,v_1,u_2,v_2$ is the solution for $c=2$. Denoting by $u_{h+1}=\big(B^y\big)^{-c}x\big(B^y\big)^c$, $v_{h+1}=\big(B^y\big)^{-c}y^2\big(B^y\big)^c$ we have
$[u_{h+1},v_{h+1}]=\big(B^y\big)^{-c}[x, y^2]\big(B^y\big)^c$ 
and therefore 
$$[u_1,v_1]\dots[u_{h+1},v_{h+1}]=\big(B\big)^c\big(B^y\big)^c\big(B^y\big)^{-c}[x, y^2]\big(B^y\big)^c=\big(B\big)^{c+1}\big(B^y\big)^{c+1}.$$
 This solution is obviously primitive since $p(u_1), p(v_1), p(u_2), p(v_2)$ generate $F_2/[F_2,F_2]$. The basis is proved. For the induction step we consider two similar cases depending on the parity of~$l$.

\textit{Case 1: $l=2n$ is even.} In this case $h=n(c-1)+1$. We will construct a primitive solution which satisfies the condition $u_1=x$, $v_1=y^{l+1}$. If $l=2$, then the statement follows from the basis of induction. By the induction hypothesis we  have a primitive solution $a_1=x$, $b_1=y^{l+1},a_2,b_2,\dots,a_{h_1},b_{h_1}$ of equation (\ref{eq1}) for $l=2n$, $h_1=n(c-1)+1$. Also by induction hypothesis we have a primitive solution $r_1=x,s_1=y^2, r_2, s_2, \dots, r_{h_2}, v_{h_2}$ of equation (\ref{eq1}) for $l=2$, $h_2=c$. If we denote by 
\begin{align}
\notag u_1&=x&&\\
\notag v_1&=y^{2n+3}&&\\
\notag u_{j}&=\big(B^{y^{2n+2}}\big)^{-1}\big(B^{y^{2n+1}}\big)^{-1}a_j\big(B^{y^{2n+1}}\big)\big(B^{y^{2n+2}}\big)&&j=2,\dots, h_1\\
\notag v_{j}&=\big(B^{y^{2n+2}}\big)^{-1}\big(B^{y^{2n+1}}\big)^{-1}b_j\big(B^{y^{2n+1}}\big)\big(B^{y^{2n+2}}\big)&&j=2,\dots, h_1\\
\notag u_{h_1+j}&=y^{2n+1}r_{j+1}y^{-(2n+1)}&&j=1,\dots, h_2-1\\
\label{vvaarr1} v_{h_1+j}&=y^{2n+1}s_{j+1}y^{-(2n+1)}&&j=1,\dots, h_2-1
\end{align}
and by $h=h_1+h_2-1=n(c-1)+1+c-1=(n+1)(c-1)+1$, then we have
\begin{align}
\notag [u_1,v_1]\dots[u_h,v_h]&=\Big([x,y^{2n+3}]\big(B^{y^{2n+2}}\big)^{-1}\big(B^{y^{2n+1}}\big)^{-1}[a_2,b_2]\dots[a_{h_1}, b_{h_1}]\Big)\\
\notag&\cdot \Big(\big(B^{y^{2n+1}}\big)\big(B^{y^{2n+2}}\big) y^{2n+1}[r_2,s_2]\dots[r_{h_2},s_{h_2}]y^{-(2n+1)}\Big)\\
\notag&=\Big([x,y^{2n+1}][a_2,b_2]\dots[a_{h_1}, b_{h_1}]\Big)\Big(y^{2n+1}BB^y[r_2,s_2]\dots[r_{h_2},s_{h_2}]y^{-(2n+1)}\Big)\\
\notag&=\Big(\big(B\big)^c\big(B^y\big)^c\dots\big(B^{y^{2n}}\big)^c\Big)\Big(y^{2n+1}\big(B\big)^c\big(B^y\big)^cy^{-(2n+1)}\Big)\\
\notag&=\big(B\big)^c\big(B^y\big)^c\dots\big(B^{y^{2n+2}}\big)^c.
\end{align}
and the statement is proved for $l=2n+2$.
The solution provided in (\ref{vvaarr1}) is primitive since the subgroup $p\Big(\big\langle u_1,v_1,\dots,u_h,v_h\big\rangle\Big)$ contains $p(u_1)=p(x)$ and $p(u_{h_1+1})=p(r_2)=p(y)^{-1}$.

{\it Case 2: $l=2n-1$ is odd.} In this case $h=\lceil (2n-1)(c-1)/2+1\rceil$. Similarly to the first case, we will show that there exists a primitive solution such that $u_1=x$, $v_1=y^l$. For $l=1$ the result follows from Corollary~\ref{main2} and Remark \ref{rmain2}. By induction hypothesis we can suppose that we have a primitive solution $a_1=x, b_1=y^l,a_2,b_2,\dots,a_{h_1},b_{h_1}$  of equation (\ref{eq1}) for  $l=2n-1$ and $h_1=\lceil (2n-1)(c-1)/2+1\rceil$. Also by induction hypothesis we can suppose that we have a solution $r_1=x,s_1=y^2, r_2, s_2, \dots, r_{h_2}, v_{h_2}$  of equation (\ref{eq1}) for $l=2$ and $h_2=c$. If we denote by 
\begin{align}
\notag u_1&=x&&\\
\notag v_1&=y^{2n+1}&&\\
\notag u_{j}&=\big(B^{y^{2n}}\big)^{-1}\big(B^{y^{2n-1}}\big)^{-1}a_j\big(B^{y^{2n-1}}\big)\big(B^{y^{2n}}\big)&&j=2,\dots, h_1\\
\notag v_{j}&=\big(B^{y^{2n}}\big)^{-1}\big(B^{y^{2n-1}}\big)^{-1}b_j\big(B^{y^{2n-1}}\big)\big(B^{y^{2n}}\big)&&j=2,\dots, h_1\\
\notag u_{h_1+j}&=y^{2n-1}r_{j+1}y^{1-2n}&&j=1,\dots, h_2-1\\
\label{vvaarr2} v_{h_1+j}&=y^{2n-1}s_{j+1}y^{1-2n}&&j=1,\dots, h_2-1
\end{align}
and by $h=h_1+h_2-1=\lceil (2n-1)(c-1)/2+1\rceil+c-1=\lceil (2n+1)(c-1)/2+1\rceil$, then we have
\begin{align}
\notag [u_1,v_1]\dots[u_h,v_h]&=\Big([x,y^{2n+1}]\big(B^{y^{2n}}\big)^{-1}\big(B^{y^{2n-1}}\big)^{-1}[a_2,b_2]\dots[a_{h_1}, b_{h_1}]\Big)\\
\notag&\cdot \Big(\big(B^{y^{2n-1}}\big)\big(B^{y^{2n}}\big) y^{2n-1}[r_2,s_2]\dots[r_{h_2},s_{h_2}]y^{1-2n}\Big)\\
\notag&=\Big([x,y^{2n-1}][a_2,b_2]\dots[a_{h_1}, b_{h_1}]\Big)\Big(y^{2n-1}BB^y[r_2,s_2]\dots[r_{h_2},s_{h_2}]y^{1-2n}\Big)\\
\notag&=\Big(\big(B\big)^c\big(B^y\big)^c\dots\big(B^{y^{2n-2}}\big)^c\Big)\Big(y^{2n-1}\big(B\big)^c\big(B^y\big)^cy^{1-2n}\Big)\\
\notag&=\big(B\big)^c\big(B^y\big)^c\dots\big(B^{y^{2n}}\big)^c
\end{align}
and the statement is proved for $l=2n+1$. The solution provided in (\ref{vvaarr2}) is primitive since the subgroup $p\Big(\big\langle u_1,v_1,\dots,u_h,v_h\big\rangle\Big)$ contains $p(u_1)=p(x)$ and $p(u_{h_1+1})=p(r_2)=p(y)^{-1}$.\hfill$\Box$
\begin{remark}{\rm 
If $c=1$, then by Proposition \ref{always} the  result of Theorem \ref{main1} holds only for $l=1$.}
\end{remark}
\begin{remark}{\rm An old problem in the geometric group theory is the problem of determining  the genus and the number $f(g)$ of Nielsen classes  for a given element $g\in [F_n,F_n]$ (see  \cite[Section~3.3]{Bes} for the definition of Nielsen classes and \cite[Question 3.11]{Bes} for the related question). If $n=2$, then Corollary \ref{main2} and Theorem \ref{main1} guarantee that for $g=([x,y])^c ({[x,y]}^y)^c\dots ({[x,y]}^{y^{l-1}})^c$ the number of Nielsen classes  is at  least $2$.}
\end{remark}
Using equality (\ref{wecken}) and construction of the map $f$ (obtained from the solution) described in Section \ref{sec1} and in details in \cite[Proposition 4.2]{GonZie} we have the following corollary.  
\begin{corollary}\label{main01} Let $c>1$, $l\geq1$ be integers and $h$ be a minimal number which satisfies the inequality $h \geq l(c-1)/2+1$. Then there exists a map $f: S_h \to T$ with $A(f)=lc$,  $MR[f]=l$, $NR[f]=1$ and  the only  root class has index $lc$. So, the Wecken property does not hold for $f$.
\end{corollary}
\section{The right part has the form $\Big([x,y]\Big)^{k+l} \Big({[x,y]}^{y}\Big)^{k-l}$}\label{var2}
The purpose of this section is to give an explicit solution for  the equation
\begin{equation}\label{eq2}
[u_1,v_1]\dots[u_h,v_h] = \Big([x,y]\Big)^{k+l} \Big({[x,y]}^{y}\Big)^{k-l}
\end{equation}
for $h=k$. If $h < k$, then by Proposition \ref{nonex} equation (\ref{eq2}) does not have solutions. The main result of this section is the following theorem
\begin{theorem}\label{main3}
Let $k>l\geq0$ be integers. Then the equation
$$[u_1,v_1]...[u_h, v_h]= \Big([x,y]\Big)^{k+l} \Big({[x,y]}^{y}\Big)^{k-l}$$ 
with unknowns $u_1,v_1,\dots,u_h, v_h$ for $h=k$ has an explicit primitive solution. Moreover, if $l\neq 0$, then every solution of this equation is primitive.
\end{theorem}
\textbf{Proof.} At first, we will prove the moreover part of the theorem. Let $u_1,v_1,\dots,u_k,v_k$ be a  solution of equation (\ref{eq2}) for $l\neq0$ and let $H=\langle u_1,v_1,\dots,u_k,v_k\rangle$. By Proposition \ref{crit}, the index $|p(F_2):p(H)|$ is equal to $1$ or $2$, and we need to prove that this index is equal to $1$. By contrary, suppose that  $|p(F_2):p(H)|=2$. From equality (\ref{wecken}) follows that the map $f:S_k\to T$ (obtained from the solution $u_1,v_1,\dots,u_k,v_k$) described in Section \ref{sec1} has two essential Nielsen root classes, one of this classes has the index $k+l$ and another one has the index $k-l$. If $l\neq 0$, then $k-l\neq k+l$, but the indices of all essential Nielsen root classes must coincide . We have a contradiction.

In order to introduce the solution of equation (\ref{eq2}) denote by
\begin{align}
\notag r_{i}&=yxy^{i-1}x^{-1}y^{-1}xy^{-i+1}x^{-1}y^{-1}&&i=1,\dots,l\\
\notag s_{i} &= yx y^{i-1}x^{-1}B^{l-i+1}yB^{k-l}y^{-i}x^2y^{-i+1}x^{-1}y^{-1}&&i=1,\dots,l\\
\notag r_{l+j}&=yxy^{(l+j)}x^{-1}y^{-1}xy^{(-l-j)}x^{-1}y^{-1}&&j=1,\dots,k-l-1\\
\notag s_{l+j} &= yxy^{(l+j)}x^{-1}B^{k-l-j}y^{-l-j}x^2y^{(-l-j)}x^{-1}y^{-1}&&j=1,\dots,k-l-1\\
\notag r_{k}&=yxy^{(k+1)}x^{-1}y^{-1}&&\\
\label{for3} s_{k}&=yxy^{-1}x^{-1}yx^{-1}y^{-1}&&
\end{align}
and let us, at first, prove some auxiliary equalities involving $r_1,s_1,\dots,r_k,s_k$.
Using induction on the number $t=1,\dots, l$ let us prove that 
\begin{equation}\label{p2}
[r_1,s_1]\dots[r_{t},s_{t}]=B^lyB^{k-l}y^{-1}B^ty^{t+1}B^{l-k}y^{-1}B^{t-1-l}xy^{-t+1}x^{-1}y^{-1}
\end{equation} 
The basis of induction ($t=1$) is proved in the following equality
\begin{align}
\notag[r_1,s_1]&=[y^{-1},yB^{l}yB^{k-l}y^{-1}xy^{-1}]\\
\notag&=y^{-1}yB^{l}yB^{k-l}y^{-1}xy^{-1}yyx^{-1}yB^{l-k}y^{-1}B^{-l}y^{-1}\\
\notag&=B^{l}yB^{k-l}y^{-1}By^2B^{l-k}y^{-1}B^{-l}y^{-1}
\end{align}
The step of induction (omitting some detailed calculations) follows from the following equality
\begin{align}
\notag[r_1,s_1]\dots[r_{t+1},s_{t+1}]&=\Big([r_1,s_1]\dots[r_{t},s_{t}]\Big)[r_{t+1},s_{t+1}]\\
\notag&=B^lyB^{k-l}y^{-1}B^ty^{t+1}B^{l-k}y^{-1}B^{t-1-l}xy^{-t+1}x^{-1}y^{-1}\\
\notag&\cdot yxy^tx^{-1}\Big[y^{-1}, B^{l-t}yB^{k-l}y^{-t-1}x\Big]xy^{-t}x^{-1}y^{-1}\\
\notag&=B^lyB^{k-l}y^{-1}B^{t+1}y^{t+2}B^{l-k}y^{-1}B^{t-l}xy^{-t}x^{-1}y^{-1}
\end{align}

Similarly to equation (\ref{p2}) using induction on the number $t=1,\dots, k-l-1$ we can prove the following equality.
\begin{equation}\label{p3}
[r_{l+1},s_{l+1}]\dots[r_{l+t},s_{l+t}]=yxy^{l+1}x^{-1}\Big(y^{-1}B^{k-l-1}y^{-l-1}B^ty^{l+t+1}B^{l+t-k}\Big)xy^{-l-t}x^{-1}y^{-1}
\end{equation} 
We will not show the proof of (\ref{p3}) here since it repeats the proof of (\ref{p2}) almost completely.

Multiplying equality (\ref{p2}) for $t=l$, equality (\ref{p3}) for $t=k-l-1$ and the value $[r_k,s_k]$ using formula (\ref{for3}) after some simple calculations we conclude that $[r_1,s_1]\dots[r_k,s_k]=B^l\big(B^y\big)^{k-l}B^k$. From this equality follows that if for $i=1,\dots,k$ we denote by 
\begin{equation}\label{for4}
u_u=B^{k}r_iB^{-k}, v_i=B^{k}s_iB^{-k},
\end{equation}
then $[u_1,v_1]\dots[u_k,v_k]=B^{k+l}\big(B^y\big)^{k-l}$, i.~e. $u_1,v_1,\dots,u_k,v_k$ is the solution of equation (\ref{eq2}). 
The images of elements $u_1,v_1,\dots,u_k,v_k$ under the homomorphism $p:F_2\to F_2/[F_2,F_2]$ generate $p(F_2)$ since $p(u_1)=p(y)^{-1}$, $p(v_1)=p(x)$. Therefore $u_1,v_1,\dots,u_k,v_k$ is a primitive solution of equation (\ref{eq2}).
 \hfill $\Box$
\begin{remark}{\rm The same result for $k=l$ follows from Theorem \ref{main1}.}
\end{remark}

~\\
Daciberg Lima Gon\c{c}alves\\
Department of Mathematics-IME, University of S\~{a}o Paulo\\
05508-090, Rua do 
Mat\~{a}o 1010, Butanta-S\~{a}o Paulo-SP, Brazil\\
email: \texttt{dlgoncal@ime.usp.br}\\
~\\
\noindent Timur Nasybullov\\
Department of Mathematics, KU Leuven KULAK\\
8500, Etienne Sabbelaan 53, Kortrijk, Belgium\\
email: \texttt{timur.nasybullov@mail.ru}


{\small
\begin{thebibliography}{00}
\bibitem{Bes}
M.~Bestvina, Questions in geometric group theory, 
\href{http://www.math.utah.edu/~bestvina}{http://www.math.utah.edu/$\sim$bestvina}.
\bibitem{BogGonKudZie1}
 S.~Bogatyi, D.~Gon\c{c}alves, E.~Kudryavtseva, H.~Zieschang, Realization of primitive branched coverings over surfaces
      following the Hurwitz approach,
Central Europ. J. of Mathematics,  V.~1, N.~2, 2003, 184--197.
\bibitem{BogGonKudZie2} S.~Bogatyi, D.~Gon\c{c}alves, E.~Kudryavtseva, H.~Zieschang, The minimal number of preimages under mappings of surfaces,  Math. Notes, V.~75, N.~1-2, 2004, 13--18.

\bibitem{BogGonKudZie3} S.~Bogatyi, D.~Gon\c{c}alves, E.~Kudryavtseva, H.~Zieschang, Realization of primitive branched coverings over closed surfaces,  Advances in Topological quantum field theory NATO Science Series, II Mathematics, Physics and Chemistry, V.~179, 2004, 297--316.
\bibitem{BogGonZie} S.~Bogatyi, D.~Gon\c{c}alves, H.~Zieschang, The minimal number of roots of surface mappings and quadratic equations in free groups, Math. Z., V.~236, N.~3, 2001, 419--452.

\bibitem{BrBPT}  R.~Brooks, R.~Brown,  J.~Pak, D.~Taylor, 
Nielsen numbers of maps of tori, Proc. Amer. Math. Soc.,  V.~52,   1975, 398--400. 

\bibitem{Cul}
M.~Culler, Using surfaces to solve equations in free groups, Topology, V.~20, N.~2, 1981, 133--145.

\bibitem{Gon} D.~Gon\c{c}alves, Coincidence of maps between surfaces, J.~Korean Math. Soc., V.~36, N.~2, 1999, 243--256. 
\bibitem{GonKudZie1} D.~Gon\c{c}alves, E.~Kudryavtseva, H.~Zieschang, Intersection index of curves on surfaces and applications to quadratic equations in free groups,  Atti Sem. Mat. Fis. Univ. Modena, V.~49, 2001, 339--400.
\bibitem{GonKudZie2}
D.~Gon\c{c}alves, E.~Kudryavtseva, H.~Zieschang, Roots of mappings on nonorientable surfaces and equations in free groups, Manuscripta Math., V.~107, N.~3, 2002, 311--341.
\bibitem{GonKudZie3} D.~Gon\c{c}alves,  E.~Kudryavtseva, H.~Zieschang,  Some quadratic equations in free group of rank $2$, Geom. Topol. Monogr., V.~14, 2008, 219--294.
\bibitem{GonZie} D.~Gon\c{c}alves, H.~Zieschang, Equations in
free groups and coincidence of mappings on surfaces, Math. Z., V.~237, N.~1, 2001, 1--29.
\bibitem{Hme1}
Ju.~Hmelevskii, Systems of equations in a free group I, Izv. Akad. Nauk SSSR Ser. Mat.,
V.~35, 1971, 1237--1268.
\bibitem{Hme2}
Ju.~Hmelevskii, Systems of equations in a free group II, Izv. Akad. Nauk SSSR Ser.
Mat., V.~36, 1972, 110--179.
\bibitem{Hme3}
Ju.~Hmelevskii, Equations in free semigroups, Proceedings of the Steklov Institute of Mathematics, N.~107, 1971, American Mathematical Society, Providence, R.I., 1976. 
\bibitem{Ki}   T. Kiang, The theory of fixed point classes, Springer-Verlag, 1989.
\bibitem{KudWeiZie} E.~Kudryavtseva, R.~Weidmann, H.~Zieschang, Quadratic equations in free groups and topological applications, Res. Exp. Math., V.~27, 2003, 83--122.
\bibitem{Lyn1}
R.~C.~Lyndon, The equation $a^2b^2=c^2$ in free groups, Michigan Math.~J,  V.~6, 1959, 89--95.
\bibitem{Lyn2}
R.~C.~Lyndon, 	Equations in free groups, Trans. Amer. Math. Soc., V.~96, 1960, 445--457.
\bibitem{Mak}
G.~Makanin, Equations in a free group, Izv. Akad. Nauk SSSR Ser. Mat., V.~46, N.~6, 1982,
1199--1273.
\bibitem{OsbZie}
R.~P.~Osborne, H.~Zieschang, Primitives in the free group on two generators, Invent.
Math., V.~63, N.~1, 1981, 17--24.
\bibitem{Seg}
D.~Segal, Words: notes on verbal width in groups, London Mathematical Society Lecture Note Series, V.~361. Cambridge University Press, Cambridge, 2009.
\bibitem{Ste}
A.~Steinberg, On equations in free groups, Michigan Math.~J., V.~18, 1971, 87--95.
\bibitem{Vdo} A.~Vdovina, Product of commutators in free products, Int. J. Algebra and Compt., V.~7, N.~4, 1997, 471--485.
\bibitem{Wic1} M.~J.~Wicks, Commutators in free products, J.~London Math. Soc., V.~37, 1962, 433--444.
\end{thebibliography}
}
\end{document}